\documentclass[11pt,reqno]{amsart}
\usepackage[T1]{fontenc}
\usepackage{lmodern}
\usepackage{amsmath,amssymb,amsthm,mathtools}
\usepackage{graphicx}
\usepackage{tikz}
\usepackage{needspace}
\usepackage{float}
\usepackage{enumitem}
\usepackage[a4paper,margin=2.65cm]{geometry}
\usepackage[hidelinks]{hyperref}
\allowdisplaybreaks
\newtheorem{theorem}{Theorem}[section]
\newtheorem{proposition}[theorem]{Proposition}
\newtheorem{remark}[theorem]{Remark}
\newtheorem{corollary}[theorem]{Corollary}
\newtheorem{construction}[theorem]{Construction}
\theoremstyle{remark}
\newcommand{\R}{\mathbb R}
\newcommand{\C}{\mathbb C}

\newcommand{\length}{\operatorname{length}}

\newcommand{\rank}{\operatorname{rank}}

\newcommand{\eps}{\varepsilon}

\title[A corank one LNE surface with non-trivial link]{A CORANK ONE LIPSCHITZ NORMALLY EMBEDDED\\ SURFACE IN $\R^4$ WITH NON-TRIVIAL LINK}
\author[L. Birbrair]{L. Birbrair}
\address{Departamento de Matem\'atica, Universidade Federal do Cear\'a, Campus do Pici, Bloco 914, Fortaleza, CE, CEP 60455-760, Brazil}
\email{lev.birbrair@gmail.com}

\author[Rodrigo Mendes]{Rodrigo Mendes}
\address{Instituto de Ci\^encias Exatas e da Natureza, Universidade da Integra\c{c}\~ao Internacional da Lusofonia Afro-Brasileira (Unilab), Campus dos Palmares, Acarape, CE, CEP 62785-000, Brazil}
\email{rodrigomendes@unilab.edu.br}

\author[J. J. Nu\~no-Ballesteros]{J. J. Nu\~no-Ballesteros}
\address{Departament de Matem\`atiques, Universitat de Val\`encia, Campus de Burjassot, 46100 Burjassot, Spain}
\email{Juan.Nuno@uv.es}

\subjclass[2020]{14B05, 32S50, 58K15}
\keywords{Lipschitz normal embedding, corank-one surface singularities, non-trivial knot links, finite $\mathcal A$-determinacy}
\dedicatory{Dedicated to the memory of C. T. C. Wall}
\date{}
\begin{document}
\begin{abstract}
We show that the knot-triviality theorem for Lipschitz normally embedded (LNE) corank one surface singularities in $\R^4$ proved in \cite{MNB16} does not generalize unrestrictedly. More precisely, we exhibit a polynomial map germ $F:(\R^2,0)\to(\R^4,0)$ of corank one at the origin, with isolated singularity and $2$-jet $(x,0,0,0)$, whose image $X=F(\R^2)$ is LNE at the origin, yet its link is the granny knot $T\#T$, and hence non-trivial. The example is obtained via a general uniform-degree homogenization procedure applied to the polynomial parametrization of the trefoil knot studied in \cite{MNB16}. We also give a finitely $\mathcal A$-determined polynomial realization with the same LNE and knot properties.
\end{abstract}
\maketitle

\section{Introduction}

Let $(X,0) \subset (\R^n,0)$ be a path-connected subanalytic set germ. Denote by $d_X$ the \emph{intrinsic (inner) metric} on $X$ (more precisely, on a chosen small neighborhood $U$ of the origin in $X$), defined as the infimum of lengths of rectifiable paths in $X$ (respectively in $U$) joining two points. We say that $(X,0)$ is \emph{Lipschitz normally embedded} (LNE) if there exist a neighborhood $U$ of $0$ in $X$ and a constant $C \ge 1$ such that
\[
d_X(z,w) \le C\|z - w\| \qquad \text{for all } z,w \in U,
\]
where $\|\cdot\|$ denotes the Euclidean norm in $\R^n$.

Let $F : (\R^2,0) \to (\R^4,0)$ be a real analytic map germ with image $X = F(\R^2)$. We say that $F$ has an \emph{isolated singularity} at the origin if, for a sufficiently small representative, $F$ is a topological embedding on its domain and an immersion outside the origin (in particular, $F^{-1}(0) = \{0\}$ and $\mathrm{rank}(dF_p) = 2$ for $p \neq 0$). By the Cone Structure Theorem proved in~\cite[Theorem~3.2]{MNB20}, for any real analytic map germ $F : (\R^2,0) \to (\R^4,0)$ with isolated singularity, the link $K_\varepsilon = X \cap S^3_\varepsilon$ is a smoothly embedded circle whose isotopy type in $S^3_\varepsilon$ is well-defined and independent of $\varepsilon > 0$ sufficiently small, providing a complete invariant for topological equivalence. In~\cite{MNB16}, explicit polynomial parametrizations realizing non-trivial knots—in particular the trefoil knot—were constructed, and the metric behavior of this construction was analyzed in~\cite{BMNB18}, where it was shown to be non-LNE.

In~\cite{BMNB18}, Birbrair, Mendes, and Nuño-Ballesteros proved that, for a broad class of real analytic surface germs of corank one in $\R^4$, the LNE property forces the link to be the unknot. It remained open whether this implication holds for all corank one real analytic surface germs with isolated singularity in $\R^4$.

The purpose of this note is to give a negative answer to this question:

\begin{theorem}[Main Theorem]\label{thm:main}
There exists a polynomial map $F : (\R^2,0) \to (\R^4,0)$ of rank $1$ at the origin, with $2$-jet $j^2 F(0) = (x,0,0,0)$, which is a global topological embedding of $\R^2$ into $\R^4$ and an immersion on $\R^2 \setminus \{(0,0)\}$ (in particular, an isolated singularity), whose image $X = F(\R^2)$ is LNE at the origin, while its link $K_\varepsilon = X \cap S^3_\varepsilon$ is isotopic to the granny knot $T \# T$ (non-trivial). In particular, the unknotting theorem of~\cite{BMNB18} does not generalize to arbitrary corank one singularities.
\end{theorem}

On the other hand, the paper \cite{BDMS25} is devoted to surfaces containing micro-knots, i.e., LNE surface germs $X\subset(\R^4,0)$ whose tangent cones are surfaces with trivial knot links, while the links of $X$ are non-trivial knots. The example of the paper shows  that the surfaces, containing micro-knots may be obtained by algebraic parametrizations.

\noindent\textbf{Asymptotic notation.} Given two real-valued functions $f, g$ defined near $0$ (or near $\infty$), we write $f \sim g$ if there exist positive constants $c_1, c_2 > 0$ such that $c_1 |g| \le |f| \le c_2 |g|$ in the corresponding regime. We write $f=O(g)$ if $|f/g|$ is bounded, and $f=o(g)$ if $f/g\to0$ (where $g\ne0$).

\noindent\textbf{Tangent cone.} Let $(X,0)\subset(\R^n,0)$ be a semialgebraic germ. Its tangent cone $C_0X$ is the closed cone over the tangent link
\[
 S_0X:=\left\{v\in\R^n:\ v=\lim_{j\to\infty}\frac{x_j}{\|x_j\|},\quad x_j\in X\setminus\{0\},\ x_j\to0\right\}.
\]
The tangent cone can equivalently be described by the tangent directions of semialgebraic arcs through the origin.

The paper is organized as follows. In Section~\ref{sec:sectional}, we establish a standalone sectional criterion for LNE of rank-one parametrized analytic surfaces in $\R^n$ (Theorem~\ref{thm:sectional}), showing via~\cite{MS24} that LNE is equivalent to uniform LNE of the lateral sections orthogonal to the $x$-axis, while the complementary horizontal sections need not be checked. In Section~\ref{sec:homog}, we introduce a uniform-degree homogenization procedure $\mathcal{H}_d(\gamma)$ for polynomial long knots $\gamma : \R \to \R^3$ and prove (Theorem~\ref{thm:construction}) that it always yields an LNE surface germ whose link is $K \# K$. 
Finally, in Section ~\ref{sec:exemplo}, we apply this construction to the classical polynomial trefoil parametrization from~\cite{MNB16,BMNB18}, proving Theorem~\ref{thm:main}. We also exhibit a finitely
$\mathcal A$-determined realization
(Theorem~\ref{thm:finite}) and show that both examples contain microknots (Corollary~\ref{cor:microknot}).

\section{Sectional criteria for Lipschitz normal embedding}\label{sec:sectional}
We begin by recalling the link criterion of Mendes and Sampaio \cite{MS24} and its formulation in the $\ell^\infty$-norm.
\begin{proposition}[Mendes--Sampaio \cite{MS24}]\label{prop:MS}
Let $(X,0)\subset(\R^n,0)$ be a closed subanalytic germ such that $(X\setminus\{0\},0)$ is connected. For any subanalytic norm $\rho$ on $\R^n$, the following are equivalent:
\begin{enumerate}[label=(\arabic*),leftmargin=1.1cm,nosep]
\item $(X,0)$ is Lipschitz normally embedded (LNE);
\item there exist $C\ge1$ and $\eps_0>0$ such that the family of links $\{K^\rho_\eps\}_{0<\eps<\eps_0}$, where $K^\rho_\eps:=X\cap\{\rho=\eps\}$, is $C$-LNE.
\end{enumerate}
Here $C$-LNE means that \emph{every member of the family is LNE with the same constant $C$}. In particular, one may take $\rho(z)=\|z\|_\infty:=\max_{1\le j\le n}|z_j|$.
\end{proposition}

\begin{remark}[Hypercube sections]\label{rem:cube}
For $\rho=\|\cdot\|_\infty$, write
\[
 Q_\eps=[-\eps,\eps]^n,\quad K^\infty_\eps=X\cap\partial Q_\eps,\quad
 A^\sigma_{j,\eps}:=X\cap Q_\eps\cap\{z_j=\sigma\eps\},
 \quad j=1,\dots,n,\ \sigma\in\{+,-\}.
\]
Let $J$ consist of those pairs $(j,\sigma)$ for which the face $\{v:\|v\|_\infty=1,\ v_j=\sigma\}$ intersects $C_0X$. By compactness of the normalized directions, for sufficiently small $\eps$ the only potentially nonempty face sections are $A^\sigma_{j,\eps}$ with $(j,\sigma)\in J$. Therefore Proposition~\ref{prop:MS} may be written as
\[
 (X,0)\text{ is LNE}\quad\Longleftrightarrow\quad
 \left\{\;\bigcup_{(j,\sigma)\in J} A^\sigma_{j,\eps}\;\right\}_{0<\eps<\eps_0}
 \text{ is }C\text{-LNE for some }C\ge1.
\]
The uniform LNE condition here applies to the \emph{union} of the relevant faces; it cannot in general be replaced by separate LNE conditions on the faces without a gluing argument.
\end{remark}

We now give the genuinely sectional reduction. All sections below are truncated by the same hypercube $Q_\eps$; this convention is essential. We write $dP_0=0$ for the gradient of $P$ and $dG_0=0$ for the Jacobian matrix of $G$ at the origin.

\begin{theorem}[Rank-one sectional LNE criterion]\label{thm:sectional}
Let $F:(\R^2,0)\to(\R^n,0)$ be real analytic, a homeomorphism onto its image germ $X=F(\R^2)$, with $\rank(dF_0)=1$. After analytic source coordinates and a suitable analytic change of ambient coordinates, write
\[
 F(x,y)=\bigl(x,P(x,y),G(x,y)\bigr),\qquad
 dP_0=dG_0=0,\qquad P(x,0)=0,\quad G(x,0)=(0,\ldots,0),
\]
where $G$ has $n-2$ coordinates, and set
\[
 \Gamma^\pm_\eps:=X\cap Q_\eps\cap\{x=\pm\eps\}.
\]
Then the following are equivalent:
\begin{enumerate}[label=(\roman*),leftmargin=1.05cm,nosep]
\item $(X,0)$ is LNE;
\item the family of hypercube links $\{K^\infty_\eps\}_{0<\eps<\eps_0}$ is $C$-LNE for some $C$;
\item the two families of lateral sections $\{\Gamma^+_\eps\}_{0<\eps<\eps_0}$ and $\{\Gamma^-_\eps\}_{0<\eps<\eps_0}$ are $C$-LNE for some common constant $C$.
\end{enumerate}
In particular, the horizontal sections orthogonal to the $P$-axis do not need to be checked to verify LNE. No immersion hypothesis outside the origin is needed for this criterion.
\end{theorem}

\begin{proof}
Since $F$ is a homeomorphism onto its image, $F(0,y)$ is not constant. Choose $P$ of minimal transverse order $p\ge2$ along $x=0$ and, if necessary, subtract its leading $y^p$ term from the components of $G$ and reverse its sign. Thus
\[
 P(0,y)=ay^p+O(y^{p+1}),\qquad G(0,y)=O(y^{p+1}),\qquad a>0.
\]
With these normalizations, analyticity gives
\begin{equation}\label{eq:taylors}
 \begin{aligned}
 P(x,y)&=ay^p+O\bigl(|y|^{p+1}+|x|\,|y|\bigr),\\
 G(x,y)&=O\bigl(|y|^{p+1}+|x|\,|y|\bigr).
 \end{aligned}
\end{equation}
If $p$ is even, the tangent cone is a half-plane, and neither (i) nor (iii) holds (the even-order obstruction of \cite[Corollary 3.4]{BMNB18}, also visible on the lateral sections). Hence we may assume $p$ is odd; in this case $C_0X=\R^2\times\{0\}^{n-2}$. By Remark~\ref{rem:cube}, only the faces $x=\pm\eps$ and $P=\pm\eps$ contribute to the hypercube link.

Put $H^\pm_\eps:=X\cap Q_\eps\cap\{P=\pm\eps\}$. We first show that these horizontal caps are uniformly Lipschitz graphs. For each $\sigma\in\{+,-\}$ and $|x|\le\eps$, the equation $P(x,y)=\sigma\eps$ in a small representative forces
\begin{equation}\label{eq:capscale}
 |y|^p\sim\eps,\qquad |x|\le\eps\lesssim|y|^p.
\end{equation}
Indeed, by \eqref{eq:taylors}, the error is $o(\eps+|y|^p)$ in this regime. Differentiating gives
\begin{equation}\label{eq:py}
 P_y(x,y)=ap y^{p-1}+O(|y|^p+|x|)\sim |y|^{p-1}>0
 \quad\bigl(P(x,y)=\sigma\eps,\ |x|\le\eps\bigr).
\end{equation}
Since $p-1$ is even, there is exactly one solution $y=\psi^\sigma_\eps(x)$ for each $|x|\le\eps$; the Implicit Function Theorem makes $\psi^\sigma_\eps$ analytic. Moreover,
\[
 |P(x,y)|\le\eps\quad\Longleftrightarrow\quad
 \psi^-_\eps(x)\le y\le\psi^+_\eps(x).
\]
Thus the hypercube link is a simple closed curve made up of the two lateral sections and the two horizontal caps.

The estimates in \eqref{eq:taylors} also give
\begin{equation}\label{eq:derivatives}
 |P_x|+\|G_x\|\le C|y|,\qquad
 \|G_y\|\le C(|x|+|y|^p).
\end{equation}
Since $(\psi^\sigma_\eps)'=-P_x/P_y$, the chain rule, \eqref{eq:capscale}--\eqref{eq:derivatives}, and $|x|\lesssim|y|^p$ give, along each cap,
\[
 \left\|\frac{d}{dx}G\bigl(x,\psi^\sigma_\eps(x)\bigr)\right\|
 \le\|G_x\|+\|G_y\|\frac{|P_x|}{|P_y|}
 \le C_1(|y|+|y|^2)\le M,
\]
where $M$ is independent of $\eps$ and $\sigma$. We can now parametrize all four arcs explicitly:
\Needspace{9\baselineskip}
\begin{align*}
 \Gamma^\sigma_\eps&:\quad y\longmapsto
 \ell^\sigma_\eps(y):=\bigl(\sigma\eps,P(\sigma\eps,y),G(\sigma\eps,y)\bigr),
 &&\psi^-_\eps(\sigma\eps)\le y\le\psi^+_\eps(\sigma\eps),\\
 H^\sigma_\eps&:\quad x\longmapsto
 h^\sigma_\eps(x):=\bigl(x,\sigma\eps,G(x,\psi^\sigma_\eps(x))\bigr),
 &&-\eps\le x\le\eps,
\end{align*}
for $\sigma\in\{+,-\}$. In particular, the caps $H^\sigma_\eps$ are uniformly Lipschitz graphs, and $K^\infty_\eps=\Gamma^+_\eps\cup\Gamma^-_\eps\cup H^+_\eps\cup H^-_\eps$.
For $z=h^\sigma_\eps(x_1)$ and $w=h^\sigma_\eps(x_2)$, we have
\begin{equation}\label{eq:graphlength}
 \length([z,w]_{H^\sigma_\eps})
 \le\sqrt{1+M^2}\,|x_1-x_2|
 \le\sqrt{1+M^2}\,\|z-w\|.
\end{equation}

\begin{figure}[htbp]
\centering
\begin{tikzpicture}[scale=0.95,>=stealth]
 \fill[black!3] (-1.8,-1.65) rectangle (1.8,1.65);
 \draw[->,gray] (-2.5,0)--(2.55,0) node[right] {$x$};
 \draw[->,gray] (0,-2.05)--(0,2.13) node[above] {$P$};
 \draw[line width=1.8pt,blue!65!black] (-1.8,-1.65)--(-1.8,1.65);
 \draw[line width=1.8pt,blue!65!black] (1.8,-1.65)--(1.8,1.65);
 \draw[line width=1.8pt,black!65] (-1.8,1.65)--(1.8,1.65);
 \draw[line width=1.8pt,black!65] (-1.8,-1.65)--(1.8,-1.65);
 \foreach \x in {-1.8,1.8} {\foreach \y in {-1.65,1.65} {\fill (\x,\y) circle (1.7pt);}}
 \node[left,blue!65!black] at (-1.84,0.85) {$\Gamma^-_\eps$};
 \node[right,blue!65!black] at (1.84,0.85) {$\Gamma^+_\eps$};
 \node[above] at (0.92,1.65) {$H^+_\eps$};
 \node[below] at (0.92,-1.65) {$H^-_\eps$};
 \node[below] at (-1.8,-1.96) {$x=-\eps$};
 \node[below] at (1.8,-1.96) {$x=\eps$};
 \node[left] at (-2.05,1.65) {$P=\eps$};
 \node[left] at (-2.05,-1.65) {$P=-\eps$};
\end{tikzpicture}
\caption{Schematic projection of the four arcs of $K^\infty_\eps$ onto the $(x,P)$-plane. The two lateral sections $\Gamma^\pm_\eps$ and the two horizontal caps $H^\pm_\eps$ meet at four corners. The transverse $G$-coordinates are suppressed.}
\label{fig:sections}
\end{figure}
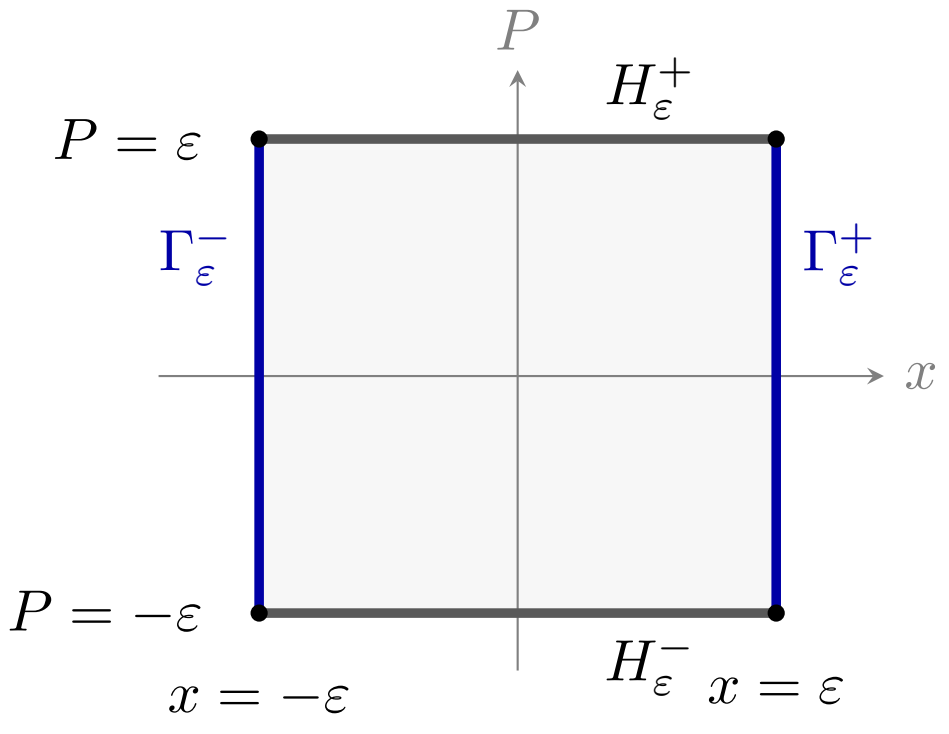

The equivalence (i)$\Leftrightarrow$(ii) is Proposition~\ref{prop:MS}. To prove (iii)$\Rightarrow$(ii), note first that $\length(\Gamma^\pm_\eps)\le A\eps$: indeed, the lateral parameter intervals lie in $|y|\le C\eps^{1/p}$ and, on $x=\pm\eps$, the bound $|P_y|+\|G_y\|\le C'(|y|^{p-1}+\eps)$ integrates to $O(\eps)$.

Let $C$ be the common LNE constant of the lateral arcs, and put $L=\sqrt{1+M^2}$. For $z\in\Gamma^\sigma_\eps$ and $w\in H^\tau_\eps$ on adjacent faces, let $c$ be the common corner. Since $x(c)=x(z)$ and the cap is an $L$-Lipschitz graph over $x$, we have $\|w-c\|\le\length([w,c])\le L|x(w)-x(c)|\le L\|z-w\|$. Using the $C$-LNE property of the lateral arc and the triangle inequality gives $d_{K^\infty_\eps}(z,w)\le C\|z-c\|+L\|z-w\|\le [C(1+L)+L]\|z-w\|$.
For opposite faces, $\|z-w\|\ge2\eps$ while a connecting path has length $O(\eps)$. Points on the same face are covered by (iii) or \eqref{eq:graphlength}. Thus the full hypercube link is uniformly LNE, proving (ii).

Conversely, if $K^\infty_\eps$ is $L$-LNE, each lateral arc is uniformly LNE: its complementary route around the link has length at least $4\eps$, whereas its own length is at most $A\eps$. Splitting according to $L\|z-w\|<4\eps$ or $L\|z-w\|\ge4\eps$ gives the common bound $\max\{L,AL/4\}\|z-w\|$ for either lateral arc. This proves (iii).
\end{proof}

\section{Uniform homogenization and the polynomial construction}\label{sec:homog}
Recall that a \emph{long knot} in $\R^3$ is a proper smooth embedding $\gamma:\R\to\R^3$ (an open, non-closed knotted curve) whose one-point compactification $\gamma(\R)\cup\{\infty\}\subset\R^3\cup\{\infty\}\cong S^3$ defines a closed knot $K\subset S^3$.

\begin{construction}[Uniform homogenization]\label{cons:homog}
Let $\gamma(u)=(p(u),q(u),r(u)):\R\to\R^3$ be an injective polynomial immersion (a polynomial long knot $K$) such that
\[
 d:=\deg r>\max\{\deg p,\deg q\},\qquad d\ge3\ \text{odd},
\]
and the leading coefficient $a$ of $r$ is positive. Define the uniform homogenization $F=H_d(\gamma):(\R^2,0)\to(\R^4,0)$ by
\[
 F(x,y):=\bigl(x,x^d r(y/x),x^d p(y/x),x^d q(y/x)\bigr)
 =:\bigl(x,P_1(x,y),P_2(x,y),P_3(x,y)\bigr),
\]
where the last three coordinates are understood as homogeneous polynomials of degree $d$ in $(x,y)$.
\end{construction}

Because $d>\max\{\deg p,\deg q\}$, we have $P_1(0,y)=ay^d$ ($a>0$), while $P_2$ and $P_3$ vanish identically on $\{x=0\}$. Since $d\ge3$, it follows that $j^2F(0)=(x,0,0,0)$ and $\rank(dF_0)=1$. Set $\widetilde\gamma:=(r,p,q)$, obtained from $\gamma$ by a cyclic permutation of its coordinates.

\begin{remark}[Hypercube link versus spherical link]\label{rem:radial}
Although the radial bi-Lipschitz homeomorphism $\Phi(z)=(\|z\|/\|z\|_\infty)z$ maps $S^3_\eps$ onto $\partial Q_\eps$, it does not preserve $X$, because $X$ is weighted homogeneous rather than conical. Nevertheless, $K_\eps=X\cap S^3_\eps$ and $K^\infty_\eps=X\cap\partial Q_\eps$ have the same knot type. For $t>0$, the scaling
\[
 \Lambda_t(x,z_1,z_2,z_3):=(tx,t^dz_1,t^dz_2,t^dz_3)
\]
satisfies $F(tx,ty)=\Lambda_t(F(x,y))$ and preserves $X$. For each $s\in[0,1]$ and $z\in X\setminus\{0\}$, the quantity $\|\Lambda_t(z)\|_s$, where $\|z\|_s=(1-s)\|z\|+s\|z\|_\infty$, increases strictly from $0$ to $+\infty$ with $t$. Let $t_s(z)>0$ be the unique continuous parameter for which $\|\Lambda_{t_s(z)}(z)\|_s=\eps$, and let $\rho(w)=\eps w/\|w\|$. Then
\[
 z\in K_\eps\longmapsto \rho\bigl(\Lambda_{t_s(z)}(z)\bigr)\in S^3_\eps
\]
is an isotopy from $K_\eps$ to the radial image $\rho(K^\infty_\eps)$.
\end{remark}

\begin{theorem}[Polynomial construction and link]\label{thm:construction}
Under the hypotheses of Construction~\ref{cons:homog}, the polynomial map $F=H_d(\gamma)$ satisfies:
\begin{enumerate}[label=(\alph*),leftmargin=1.05cm,nosep]
\item $F$ is a global topological embedding of $\R^2$ into $\R^4$ and an immersion on $\R^2\setminus\{(0,0)\}$ (in particular, $F$ has an isolated singularity at the origin);
\item $(X,0)$ is Lipschitz normally embedded, where $X=F(\R^2)$;
\item there exists $\eps_0>0$ such that, for every $0<\eps\le\eps_0$, the link $K_\eps=X\cap S^3_\eps$ is well-defined, its isotopy type is independent of $\eps$, and $K_\eps$ is isotopic to the connected sum $K\#K$ (in particular, $K_\eps$ is non-trivial whenever $K$ is non-trivial).
\end{enumerate}
\end{theorem}

\begin{proof}
\emph{(a) Isolated singularity and topological embedding.} Since $F(0,y)=(0,ay^d,0,0)$, the oddness of $d$ and the injectivity of $\gamma$ imply that $F$ is injective on $\R^2$. It is also proper, by dominance of $ay^d$ in the second coordinate when $x$ stays bounded and $|y|\to\infty$ (while $|x|\le\|F(x,y)\|$). Moreover, $\rank(dF_0)=1$ and $\rank(dF_{(x,y)})=2$ for $(x,y)\ne0$: for $x\ne0$ this follows from $\widetilde\gamma'\ne0$, while at $(0,y)$, $y\ne0$, it follows from $\partial_yP_1(0,y)=day^{d-1}\ne0$. Thus $F$ is a global topological embedding and an immersion outside $0$.

\emph{(b) Proof of the LNE property.} By (a) and $P_1(0,y)=ay^d$ with $d\ge3$ odd, Theorem~\ref{thm:sectional} applies (with $P=P_1$, $G=(P_2,P_3)$, and $p=d$): it suffices to prove uniform LNE of the two lateral sections $\Gamma^\pm_\eps=X\cap Q_\eps\cap\{x=\pm\eps\}$.

For $x=\eps>0$, setting $y=\eps u$ gives
\[
 F(\eps,\eps u)=\bigl(\eps,\eps^d\widetilde\gamma(u)\bigr),\qquad
 J_\eps:=\bigl\{u\in\R:\|\widetilde\gamma(u)\|_\infty\le\eps^{1-d}\bigr\}.
\]
We first observe that $J_\eps$ is an interval for small $\eps$. Indeed, there is $R>0$ such that $r$ is strictly increasing on both tails $(-\infty,-R]$ and $[R,\infty)$ and $|p(u)|,|q(u)|\le|r(u)|$ whenever $|u|\ge R$. For $T=\eps^{1-d}$ larger than the maximum of $|p|,|q|,|r|$ on $[-R,R]$, it follows that
\[
 J_\eps=\{u\in\R:|r(u)|\le T\}=[u_-(T),u_+(T)],
 \qquad r(u_\pm(T))=\pm T.
\]
We next show that $\widetilde\gamma(\R)$ is globally LNE. Write $r(u)=au^d+s(u)$ with $\deg s<d$. For $u<v$ and $H=\max\{|u|,|v|\}\ge1$, the oddness of $d$ gives
\[
 \frac{v^d-u^d}{v-u}\ge\frac12H^{d-1},\qquad
 |s(v)-s(u)|\le AH^{d-2}(v-u)
\]
for some $A>0$. Indeed, up to symmetry take $H=v>0$ and $u=\lambda v$, $\lambda\in[-1,1)$; then $(1-\lambda^d)/(1-\lambda)\ge1/2$. Hence, for $H\ge H_0$ sufficiently large,
\[
 r(v)-r(u)\ge\frac a4H^{d-1}(v-u).
\]
As $\|\widetilde\gamma'(t)\|\le BH^{d-1}$ for $t\in[u,v]$, we obtain
\[
 \int_u^v\|\widetilde\gamma'(t)\|\,dt
 \le\frac{4B}{a}|r(v)-r(u)|
 \le\frac{4B}{a}\|\widetilde\gamma(v)-\widetilde\gamma(u)\|.
\]
For $u,v\in[-H_0,H_0]$, the same inequality, with another constant, follows from the fact that a compact connected $C^1$ embedded regular curve is LNE, since $\widetilde\gamma$ is an injective polynomial immersion. Consequently, there exists $C\ge1$ such that
\begin{equation}\label{eq:longLNE}
 \int_u^v\|\widetilde\gamma'(t)\|\,dt
 \le C\|\widetilde\gamma(v)-\widetilde\gamma(u)\|\qquad(u<v\text{ in }\R).
\end{equation}
Since $J_\eps$ is an interval, \eqref{eq:longLNE} yields, for $u<v$ in $J_\eps$,
\[
 d_{\Gamma^+_\eps}\bigl(F(\eps,\eps u),F(\eps,\eps v)\bigr)
 \le\eps^d\int_u^v\|\widetilde\gamma'(t)\|\,dt
 \le C\|F(\eps,\eps v)-F(\eps,\eps u)\|.
\]
Finally, homogeneity and oddness of $d$ give $F(-x,-y)=-F(x,y)$; thus $\Gamma^-_\eps=-\Gamma^+_\eps$ is $C$-LNE as well. Theorem~\ref{thm:sectional} implies that $(X,0)$ is LNE.

\emph{(c) Identification of the link.} By (a) and the Cone Structure Theorem \cite[Theorem 3.2]{MNB20}, there exists $\eps_0>0$ such that the link $K_\eps=X\cap S^3_\eps$ is well-defined and has constant isotopy type for $0<\eps\le\eps_0$. By Remark~\ref{rem:radial}, it suffices to determine the knot type of $K^\infty_\eps=X\cap\partial Q_\eps$. The tangent cone is the $(x,P_1)$-plane, and for small $\eps$ the hypercube link decomposes as
\[
 K^\infty_\eps=\Gamma^+_\eps\cup\Gamma^-_\eps\cup H^+_\eps\cup H^-_\eps.
\]
Each cap $H^\sigma_\eps$ lies inside the face $\partial Q_\eps\cap\{P_1=\sigma\eps\}$ and is a graph
\[
 x\in[-\eps,\eps]\longmapsto
 \bigl(x,\sigma\eps,g_2^\sigma(x),g_3^\sigma(x)\bigr),
 \qquad |g_j^\sigma(x)|=o(\eps)<\eps.
\]
Thus each cap is an unknotted arc in its $3$-face.

The sphere $\Sigma:=\partial Q_\eps\cap\{x=0\}\cong S^2$ divides $\partial Q_\eps$ into two closed $3$-balls $B^\pm=\partial Q_\eps\cap\{\pm x\ge0\}$. The link meets $\Sigma$ transversely in the two points $q^\pm=H^\pm_\eps\cap\{x=0\}$ and splits into proper arcs $\Delta^\pm=K^\infty_\eps\cap B^\pm$. Choose a simple arc $\beta\subset\Sigma$ joining $q^+$ to $q^-$ and close the two tangles by
\[
 K_+:=\Delta^+\cup\beta,\qquad K_-:=\Delta^-\cup(-\beta).
\]
The standard connected-sum decomposition across $\Sigma$ gives $K^\infty_\eps\cong K_+\#K_-$. The antipodal map $A(z)=-z$ has $\det(-I_4)=+1$ and is orientation preserving on $\partial Q_\eps$; because $F(-x,-y)=-F(x,y)$, it interchanges the two tangles. Its image of $\beta$ is isotopic in $\Sigma$ relative to the endpoints to the reversed closing arc. Consequently $K_-\cong K_+$ as unoriented knot types, with the \emph{same chirality}.

The half-caps $H^\pm_\eps\cap\{x\ge0\}$ and $\beta$ form an unknotted closing tangle for $\Gamma^+_\eps$, disjoint from its interior. By $F(\eps,\eps u)=(\eps,\eps^d\widetilde\gamma(u))$, the lateral section is a translated and rescaled copy of $\widetilde\gamma(J_\eps)$. For $|u|$ sufficiently large, $r$ is strictly increasing and dominates $p$ and $q$ in absolute value. Hence the resulting closure is isotopic to the one-point compactification of the long knot $\gamma$. Thus $K_+\cong K$ and
\[
 K_\eps\cong K^\infty_\eps\cong K\#K.
\]
By additivity of knot genus, $g(K\#K)=2g(K)>0$ whenever $K$ is non-trivial, proving (c).
\end{proof}

\section{The explicit counterexample and microknots}\label{sec:exemplo}
\begin{proof}[Proof of Theorem~\ref{thm:main}]
Consider the classical polynomial parametrization of the long trefoil knot used in \cite{MNB16,BMNB18}:
\[
 \gamma(u)=\bigl(p(u),q(u),r(u)\bigr)
 =\left(u^3-u,\ u^4-\frac32u^2,\ u^5-\frac53u^3+\frac35u\right).
\]
We verify directly that $\gamma$ is an injective immersion on $\R$:
\begin{itemize}[nosep,leftmargin=.65cm]
\item The roots of $p'(u)=3u^2-1$ (namely $u=\pm1/\sqrt3$) do not annihilate $q'(u)=4u^3-3u$, so $\gamma'(u)\ne0$ for every $u\in\R$.
\item If $p(u_1)=p(u_2)$ and $q(u_1)=q(u_2)$ with $u_1<u_2$, then
\[
 u_1^2+u_1u_2+u_2^2=1,\qquad
 (u_2^2-u_1^2)(u_1^2+u_2^2-\tfrac32)=0.
\]
If $u_2=-u_1$, this yields $(u_1,u_2)=(-1,1)$, where $r(1)-r(-1)=-2/15\ne0$; otherwise $u_1u_2=-1/2$, $u_1+u_2=\pm1/\sqrt2$, and
\[
 r(u_2)-r(u_1)=(u_2-u_1)\left(\frac54-\frac53+\frac35\right)
 =\frac{11}{60}(u_2-u_1)\ne0.
\]
\end{itemize}
In fact, these polynomial identities hold over $\C$, so $\gamma_\C:\C\to\C^3$ is an injective immersion.

Since $d=\deg r=5>\max\{\deg p,\deg q\}=4$ is odd with leading coefficient $1>0$, applying Construction~\ref{cons:homog} produces the polynomial map
\begin{equation}\label{eq:example}
 \boxed{F(x,y)=\left(x,\ y^5-\frac53x^2y^3+\frac35x^4y,\ x^2y^3-x^4y,\ xy^4-\frac32x^3y^2\right).}
\end{equation}
By Theorem~\ref{thm:construction}, $F$ is a global topological embedding of $\R^2$ into $\R^4$, an immersion outside the origin, has $2$-jet $j^2F(0)=(x,0,0,0)$, and its image $X=F(\R^2)$ is LNE at the origin. Furthermore, the link of $X$ is the granny knot $T\#T$, which has genus $2$ and is therefore non-trivial.
\end{proof}

\begin{figure}[H]
 \centering
 \includegraphics[width=.88\linewidth]{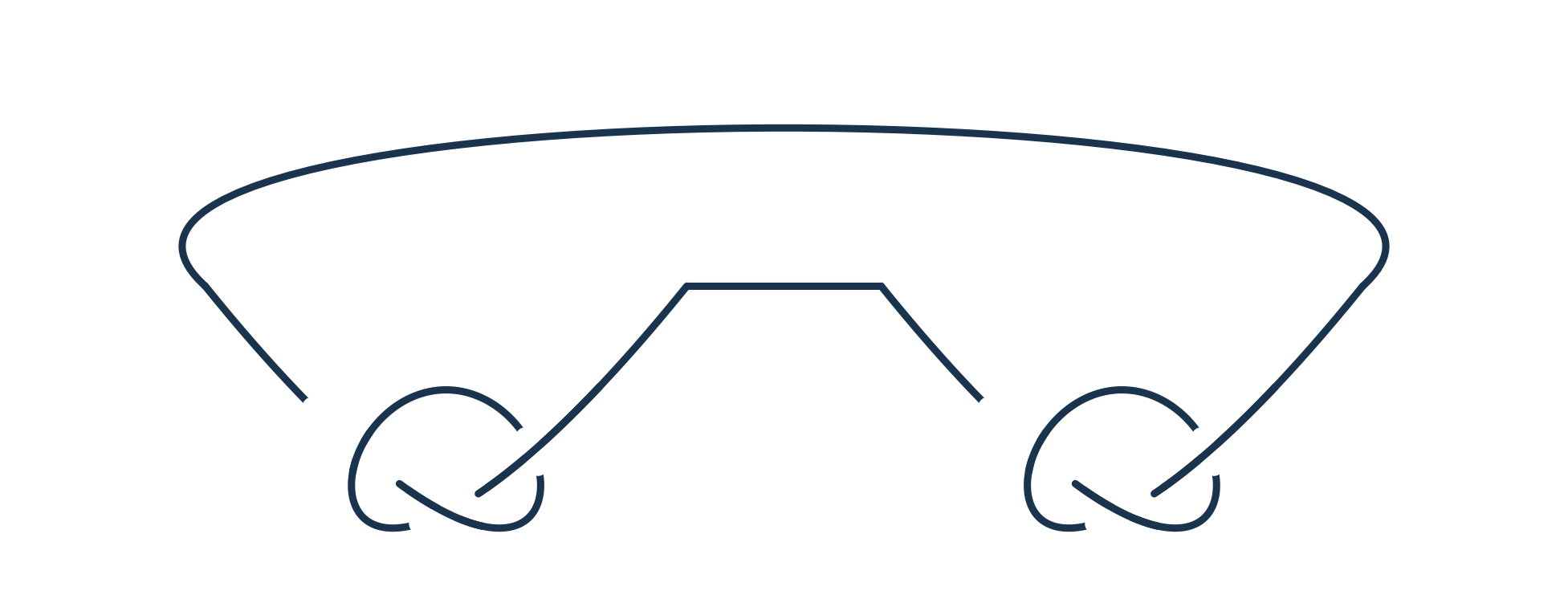}
 \caption{A diagram of the granny knot $T\#T$ (two trefoils with the same chirality). The two trefoil tangles were obtained from the planar projection $(p(u),q(u))$ of the polynomial curve $\gamma$ in \eqref{eq:example}, with the crossing information determined by $r(u)$.}
 \label{fig:granny}
\end{figure}

\begin{theorem}[A finitely determined realization]\label{thm:finite}
The polynomial germ
\begin{equation}\label{eq:finiterealization}
\widehat F(x,y)=\bigl(x,\,y^5-\tfrac53x^2y^3+\tfrac35x^4y,\,
 x^2y^3-x^4y+y^7,\,xy^4-\tfrac32x^3y^2+y^8\bigr)
\end{equation}
has an isolated singularity over $\R$, its image is LNE, and its link is the granny knot $T\#T$. Its complexification, and consequently the real germ itself, is finitely $\mathcal A$-determined.
\end{theorem}

\begin{proof}
Write $F_s=F+s(0,0,y^7,y^8)$, $0\le s\le1$, so that $F_0=F$ and $F_1=\widehat F$.
\emph{Complex finite determinacy.} The complexified polynomial curve $\widetilde\gamma_{\C}=(r,p,q)$ is an injective immersion: the algebraic identities in the proof of Theorem~\ref{thm:main} remain valid over $\C$, and $p'$ and $q'$ have no common root. Since the vector $\partial_y(P_1,P_2,P_3)$ is homogeneous of degree four and has no zero in $\C^2\setminus\{0\}$, the terms of degrees seven and eight cannot create critical points away from $0$ on a small representative.

We also exclude complex double points. Any such pair has the form $(x,y),(x,z)$ with $y\ne z$. When $x=0$, the equations $y^5=z^5$ and $y^7=z^7$ imply $y=z$. For $x\ne0$, set $u=y/x$ and $v=z/x$. If $u,v$ remain bounded, the curves $u\mapsto(r(u),p(u)+x^2u^7,q(u)+x^3u^8)$ remain embedded on compact sets, since they converge in $C^1$ to $\widetilde\gamma_\C$. Thus, in a hypothetical sequence of double points approaching $0$ not covered by this case, $|u|$ and $|v|$ become unbounded; indeed, $r(u)=r(v)$ prevents one from remaining bounded while the other diverges. After interchanging $y,z$ if necessary, put $a=x/y\to0$ and $b=z/y\to\zeta$, where equality of the second coordinates gives $\zeta^5=1$. If $\zeta=1$, the divided difference $\Delta P_1(x;y,z):=[P_1(x,y)-P_1(x,z)]/(y-z)=5y^4(1+o(1))$ cannot vanish. For $\zeta\ne1$, equality of the last two coordinates, after substituting $x=ay$, $z=by$ and dividing by $y^5$, gives
\begin{equation}\label{eq:complexdouble}
\begin{aligned}
 a^2(b^3-1)-a^4(b-1)+y^2(b^7-1)&=0,\\
 a(b^4-1)-\tfrac32a^3(b^2-1)+y^3(b^8-1)&=0.
\end{aligned}
\end{equation}
Because $\zeta^3,\zeta^4,\zeta^7,\zeta^8\ne1$, these two equations respectively yield $|a|\sim|y|$ and $|a|\sim|y|^3$, a contradiction. Thus $\widehat F_\C$ is an injective immersion on a punctured representative, hence finitely $\mathcal A$-determined by the Mather--Gaffney criterion \cite[Theorem~2.1]{Wall}. Finite determinacy of the complexification implies finite determinacy of the real analytic germ.

\emph{LNE and the real embeddings.} On the lateral sections $x=\pm\eps$, put $y=xu$. After dividing the transverse coordinates by $x^5$, we obtain $\widetilde\gamma(u)+E_{x,s}(u)$, where $E_{x,s}(u)=(0,sx^2u^7,sx^3u^8)$. Since the lateral sections lie in the hypercube $Q_\eps=[-\eps,\eps]^4$ implies $|r(u)|\le\eps^{-4}$, hence $|u|=O(\eps^{-4/5})$ and $|xu|=O(\eps^{1/5})$. Moreover, $\|E'_{x,s}(u)\|\le C(|xu|^2+|xu|^3)\|\widetilde\gamma'(u)\|\le\kappa_\eps\|\widetilde\gamma'(u)\|$, where $\kappa_\eps=O(\eps^{2/5})\to0$. For $|u|$ sufficiently large, $r$ is strictly increasing and dominates the other coordinates, including their perturbations; hence the admissible parameter set remains the interval $J_\eps$ of Theorem~\ref{thm:construction}(b). The derivative estimate preserves immersion and, together with~\eqref{eq:longLNE}, bounds the perturbed arc lengths by $C(1+\kappa_\eps)/(1-C\kappa_\eps)$ times the distance between their endpoints. Thus the lateral sections are uniformly LNE and embedded; the real maps remain embeddings on a punctured representative, and Theorem~\ref{thm:sectional} gives LNE of $F_s(\R^2)$, in particular for $s=1$.

\emph{Preservation of the knot.} For $\theta=(a,b)\in S^1$ write $F_s(t\theta)=(ta,t^5A(\theta)+st^7B_7(\theta)+st^8B_8(\theta))$. Here $A=(P_1,P_2,P_3)|_{S^1}$ is the degree-five homogeneous part, while $B_7(a,b)=(0,b^7,0)$ and $B_8(a,b)=(0,0,b^8)$ are the perturbing vectors. Since $A$ vanishes only when $b=0$, we have $\|B_7\|+\|B_8\|\le C\|A\|$. It follows, uniformly in $s$ and $\theta$, that
$\frac{d}{dt}\|F_s(t\theta)\|^2=2ta^2+10t^9\|A(\theta)\|^2+O(t^{11}\|A(\theta)\|^2)>0$ for $0<t\ll1$. Hence every source ray meets the inverse image of a fixed sufficiently small sphere exactly once. The resulting family of embedded circles in $S^3_\eps$ is an isotopy from the link of $F$ to that of $\widehat F$, so the latter is $T\#T$.
\end{proof}

We now connect this construction with the theory of microknots introduced by Birbrair, Denkowski, Medeiros, and Sampaio \cite{BDMS25}. In the present example the non-trivial tangles are located in narrowing sectors around the two opposite $x$-axis rays.
 
We say that an LNE surface germ $(X,0)$ with an isolated singularity at zero contains a microknot if it is outer bi-Lipschitz equivalent to the tangent cone of $(X,0)$ at $0\in \R^4$, but not ambient bi-Lipschitz equivalent to the tangent cone of $(X,0)$ at $0\in \R^4$.

\begin{remark}\label{contains-mikroknot-remark}\normalfont
The examples from the paper (\cite{BDMS25}) satisfied this definition. The definition is not exactly the same as in \cite{BDMS25},
but we present this definition for simplicity.
\end{remark}

\begin{corollary}[Presence of microknots]
\label{cor:microknot}
LNE surface $X$ described by
\eqref{eq:example} or \eqref{eq:finiterealization} contains microknots
corresponding to $K$ in each half-surface
$X_\pm := F(\{\pm x \ge 0\})$.
\end{corollary}

\begin{proof}
Notice that the tangent cone of $(X,0)$ is the plane $\R^2$, whose link is unknotted. Since $(X,0)$ is LNE, the theory of H\"older complexes \cite{birbrair99} gives that $(X,0)$ is outer bi-Lipschitz equivalent to $(\R^2,0)$. In both cases the link of $X$ is non-trivial, so $(X,0)$ is not ambient bi-Lipschitz equivalent to its tangent cone.
\end{proof}

\textbf{Acknowledgments}: During the preparation of this manuscript the authors used ChatGPT and Gemini to assist in refining and revising the manuscript. The authors take full responsibility for the content and accuracy of the final manuscript.

\end{document}